\documentstyle[12pt]{article}
\renewcommand{\a}{\alpha}
\renewcommand{\b}{\beta}

\newcommand{\k}{\kappa}

\renewcommand{\L}{\Lambda}

\newcommand{\th}{\theta}

\begin{document}
\setlength{\baselineskip}{18pt}

\begin{center}
{\large\bf Self-similar Intermediate
Structures}\medskip\\
{\large\bf in Turbulent Boundary Layers}\medskip\\ {\large\bf at Large
Reynolds Numbers}

\vspace{.5 truein}
{\sc G. I. Barenblatt, A. J. Chorin}

\bigskip
Department of Mathematics\\
and Lawrence Berkeley National Laboratory\\ University of California\\
Berkeley, California 94720, USA

\medskip
and \medskip

{\sc V. M. Prostokishin}

\bigskip
P.~P.~Shirshov Institute of Oceanology\\ Russian Academy of Sciences\\
36 Nakhimov Prospect\\
Moscow 117218, Russia
\end{center}

\bigskip\bigskip

{\bf Abstract.}
Processing the data from a large variety of zero-pressure-gradient boundary
layer flows shows that the Reynolds-number-dependent scaling law, which the
present authors obtained earlier for pipes, gives an accurate description
of the velocity distribution in a self-similar intermediate region of
distances from the wall adjacent to the viscous sublayer. The appropriate
length scale that enters the definition of the boundary layer Reynolds
number is found for all the flows under investigation.

Another
intermediate self-similar region between the free stream and the first
intermediate region is found under conditions of weak free stream
turbulence. The effects of turbulence in the free stream and of wall
roughness are assessed, and conclusions are drawn.

\newpage
\section{Introduction}

Asymptotic laws for
wall-bounded turbulent shear flows at large Reynolds numbers are
considered. Classical
examples of such flows are the flows in pipes, channels, and boundary
layers. This class of flows is of major fundamental and practical
importance. All these flows share as dimensional governing parameters the
shear stress at the wall $\tau$ and the fluid's properties, its density
$\rho$ and dynamic viscosity $\mu$. From these parameters two important
quantities can be formed: the {\em dynamic} or {\em friction} velocity $u_*
= (\tau/\rho)^{\frac {1}{2}}$ and the length scale $\delta = \nu/u_*$,
where $\nu = \mu/\rho$ is the fluid's kinematic viscosity. The length scale
$\delta$ is tiny at large Reynolds numbers, and in the layer where the
dimensionless distance of the wall $y/\delta$ is less than, say, 70
(viscous sublayer) the viscous stress is comparable with the Reynolds
stress created by vortices. Outside this viscous sublayer, at $y/\delta >
70$, the contribution of the viscous stress is small. We emphasize that
`small' is not always synonymous with `negligible', and indeed we will see
that here is a case where it is not.

In 1930, the great mechanician Th.~von K\'arm\'an proposed in explicit
form the hypothesis that outside the viscous sublayer the contribution of
viscosity can be neglected. On the basis of this assumption he derived the
{\it universal} (i.e.~Reynolds number independent) {\it logarithmic law}
for the distribution of the mean velocity $u$ over the cross-section:
\begin{equation}
\label{eq1}
\phi=\frac{u}{u_*}=\frac 1\k \mbox{ ln } \eta+C \ , \qquad
\eta=\frac{u_*y}{\nu}
\end{equation}
where $y$ is the distance from the wall; the constants $\k$ (the von
K\'arm\'an constant) and $C$ should be identical for all turbulent
wall-bounded shear flows at high Reynolds numbers, and the law (\ref{eq1})
should be valid in intermediate regions between, on one hand, the viscous
sublayer and, on the other, the external parts of the flows, e.g.~vicinity
of the axis in pipe flow, or vicinity of the external flow in the boundary
layer. In 1932, L.~Prandtl, the greatest mechanician of this century, came
to the law (\ref{eq1}) using a different approach, but effectively with the
same basic assumption. The law (\ref{eq1}) is known as the von
K\'arm\'an-Prandtl universal logarithmic law. More recent derivations
which, however, follow the same ideas and the same basic assumption, often
in an implicit form, can be found in monographs by Landau and Lifshits
(1987), Monin and Yaglom (1971), Schlichting (1968) and in a recent
textbook by Spurk (1997).

According to the von K\'arm\'an-Prandtl law (\ref{eq1}), all experimental
points corresponding to the intermediate region should collapse on a single
universal straight line in the traditional coordinates $\mbox{ln
}\eta,\phi$.

Subsequent investigations showed, however, that this is not what happens.
First, the experiments showed systematic deviations from the universal
logarithmic law (\ref{eq1}) even if one is willing to tolerate a variation
in the constants $\k$ and $C$ (from less than 0.4 to 0.45 for $\k$, and
from less than 5.0 to 6.3 for $C$). Furthermore, using analytic and
experimental arguments, the present authors showed [Barenblatt (1991,
1993); Barenblatt and Prostokishin (1993); Barenblatt, Chorin and
Prostokishin (1997b); Chorin (1998)] that the fundamental von K\'arm\'an
hypothesis on which the derivation of the universal law (\ref{eq1}) was
based, i.e.~the assumption that the influence of viscosity disappears
totally outside the viscous sublayer, is inadequate. In fact, this
hypothesis should be replaced by the more complicated one of incomplete
similarity, so that the influence of viscosity in the intermediate region
remains, but the viscosity enters only in power combination with other
factors. This means that the influence of the Reynolds number, i.e.~both of
the viscosity and the external length scale, e.g.~the pipe diameter,
remains and
should be taken into account in the intermediate region.

For the readers' convenience we present here briefly the concept of incomplete
similarity; a more detailed exposition can be found in Barenblatt, (1996).
The mean velocity gradient $\partial_y u$ in turbulent shear flows can be
represented in the general form suggested by dimensional analysis
\[
\partial_y u =\frac{u_*}{y} \Phi(\eta,{\mbox{Re}}) \ . \]

In the intermediate region, $\eta =u_*y/\nu$ is large, and we consider the
case of large Reynolds number. The basic von K\'arm\'an hypothesis
corresponds to the assumption that the dimensionless function
$\Phi(\eta,\mbox{Re})$ at large $\eta$ and Re can be replaced by a constant
$1/\k$, its limit as
$\eta\to\infty$, Re~$\to\infty$. This corresponds to {\it complete
similarity} both in $\eta$ and Re. The assumption of {\it incomplete
similarity} in $\eta$ means that at large $\eta$ a finite limit of the
function $\Phi$ does not exist, but that this function can be represented as
\[
\Phi =C(\mbox{Re})\eta^{\alpha(\mbox{Re})} \] i.e., the velocity gradient
has a {\it scaling} intermediate asymptotics. Here the functions $C$(Re)
and $\a$(Re) should be specified.

Using some additional analytic and experimental arguments the present
authors came to the Reynolds-number-dependent {\it scaling} law of the form
\begin{equation}
\label{eq2}
\phi=\frac{u}{u_*}=(C_0\mbox{ ln }
{\mbox{Re}}+C_1)\eta^{c/\mbox{\footnotesize ln} {\mbox{\footnotesize
\ Re}}}
\ . \end{equation}
where the constants $C_0$, $C_1$ and $\a$ must be universal. The scaling
law (\ref{eq2}) was compared with what seemed (and seems to us up to now)
to be the best available data for turbulent pipe flows, obtained by
Nikuradze (1932), under the guidance of Prandtl at his Institute in
G\"ottingen. The comparison has yielded the following values for the
coefficients \begin{equation}
\label{eq3}
c={\textstyle{\frac 32}} \ , \qquad
C_0=\frac{1}{\sqrt{3}} \ , \qquad
C_1={\textstyle{\frac 52}}
\end{equation}
when the Reynolds number Re was taken in the form \begin{equation}
\label{eq4}
{\mbox{Re}} =\frac{\bar u d}{\nu} \ .
\end{equation}
Here $\bar u$ is the average velocity (the total flux divided by the pipe
cross-section area) and $d$ is the pipe diameter. The final result has the
form
\begin{equation}
\label{eq5}
\phi =\big(\frac{1}{\sqrt{3}} \mbox{ ln} {\mbox{ Re}} + {\textstyle{\frac
52}} \big)\eta^{3/2\ \mbox{\footnotesize ln} {\mbox{\footnotesize \ Re}}}
\end{equation}
or, equivalently
\begin{equation}
\label{eq6}
\phi =\left(\frac{\sqrt{3}+5\a}{2\a}\right) \eta^{\a} \ , \qquad
\a=\frac{3}{2\mbox{ln }{\mbox{Re}}} \ .
\end{equation}
The scaling law (\ref{eq5}) produces
separate curves $\phi(\mbox{ln }\eta$,Re) in the
traditional $(\mbox{ln }\eta,\phi)$ plane, one
for each value of the Reynolds number. This is
the principal difference between the law
(\ref{eq5}) and the universal logarithmic law
(\ref{eq1}). We showed that the family
(\ref{eq5}) of curves having Re as parameter has
an envelope, and that in the $\mbox{ln}\eta,\phi$
plane this envelope is close to a straight line,
analogous to (\ref{eq1}) with the values $\k=0.4$
and $C=5.1$. Therefore, if the experimental
points are close to the envelope they can lead
to the illusion that they confirm the universal
logarithmic law (\ref{eq1}).

The Reynolds-number-dependent scaling law can be reduced to a self-similar
universal form
\begin{equation}\label{eq7}
\psi=\frac 1\a \mbox{ln }\left(
\frac{2\a\phi}{\sqrt{3}+5\a}\right)=\mbox{ln }\eta \ , \qquad
\a=\frac{3}{2\mbox{ln }{\mbox{Re}}}, \end{equation}
so that contrary to what happens in the
$(\mbox{ln }\eta,\phi)$ plane, in the
$(\mbox{ln }\eta,\psi)$ plane the experimental points should collapse onto
a single straight line---the bisectrix of the first quadrant. This
statement received a ringing confirmation from the processing of
Nikuradze's 1932 data (Barenblatt and Prostokishin (1993); Barenblatt,
Chorin and Prostokishin (1997b)).

An important remark should be made here. Izakson, Millikan and von Mises
(IMM, see, e.g.~Monin and Yaglom (1971)) gave an elegant derivation of the
universal logarithmic law based on what is now known as matched asymptotic
expansions. This derivation, which seemed to be unbreakable,
persuaded fluid dynamicists that this law was a truth which will enter
future turbulence theory essentially unchanged. In the papers of the
present authors (Barenblatt, Chorin (1996, 1997)), it was demonstrated that
the scaling law (\ref{eq2}) is compatible with the properly modified IMM
procedure. The method of vanishing viscosity (Chorin, (1988, 1994)) was
used in this
modification.

Let us turn now to shear flows other than flows in pipes. By the same
logic, the
scaling law (\ref{eq5}) should be also valid for an intermediate region
adjacent to the viscous sublayer for all good quality experiments performed
in turbulent shear flows at large Re.

The first question is, what is the appropriate definition of the Reynolds
number for these flows which will make the formula (\ref{eq5}) applicable?
This is a very important point---if the universal
Reynolds-number-independent logarithmic law were valid, the definition of
the Reynolds number would be irrelevant provided it were sufficiently
large. For the scaling law (\ref{eq5}) this is not the case. Indeed, if
the scaling law (\ref{eq5}) has general applicability it should be possible
to find, for every turbulent shear flow at large Reynolds number, an
appropriate definition of the Reynolds number which will make the scaling
law (\ref{eq5}) valid.

There exists nowadays a
large amount of data for an important class of wall-bounded turbulent shear
flows: turbulent zero-pressure-gradient boundary layers. These data were
obtained over the last 25 years by various authors using various set-ups.
For boundary layers the traditional definition of the Reynolds number is
\begin{equation}\label{eq8} {\mbox{Re}}_\th =\frac{U\th}{\nu}
\end{equation}
where $U$ is the free stream velocity, and $\th$ is a characteristic length
scale---the momentum displacement thickness. The question which we asked
ourselves was, is it possible to find for each of these flows a particular
length scale $\L$, so that the scaling law (\ref{eq5}) will be valid for
all of them with the same values of the constants. Of course, in each case
the length-scale $\L$ could be influenced by the contingencies of the
particular experiment, but the question of decisive importance is whether
such a length scale exists. The answer within the accuracy of the
experiments is affirmative.

We present here the results of the processing all the experimental data
available to us, in particular, all the data collected in a very
instructive review by Fernholz and Finley (1996). We show that for all of
these flows, without any exception, the scaling law (\ref{eq5}) is observed
with an instructive accuracy over the whole intermediate region, if the
Reynolds number is defined properly, i.e., if the characteristic length
$\L$ entering the Reynolds number
\begin{equation}\label{eq9}
{\mbox{Re}} =\frac{U\L}{\nu}
\end{equation}
is properly determined.
Moreover, we show that for all the flows where the turbulence in the
external flow is small, there exists a sharply distinguishable second
intermediate region between the first one where the scaling law (\ref{eq5})
is valid and the external homogeneous flow. The average velocity
distribution in this second intermediate region is also self-similar of
scaling type:
\begin{equation}\label{eq10}
\phi=B\eta^\b
\end{equation}
where $B$ and $\b$ are constants.

However, a Reynolds number dependence of the power $\b$ was not observed.
Within the accuracy of the experimental data $\b$ is close to 1/5. When the
turbulence in the external homogeneous flow becomes significant, the second
self-similar region deteriorates and the power $\b$ decreases with growing
external turbulence until the second intermediate region disappears
completely.

\section{The first group of zero-pressure-gradient boundary layer experiments}

We will explain later why we divided the experimental data into three
groups. Here it is sufficient to note that all available sets of
experimental data were eventually taken into account.

The original data were always presented by their authors in the form of
graphs in the traditional
$(\mbox{ln }\eta,\phi)$ plane, suggested by the universal logarithmic law
(\ref{eq1}). The shape of original graphs was always similar to the one
presented qualitatively in Figure 1a. Therefore, the first rather trivial
step was to replot the data in the doubly logarithmic coordinates
$(\mbox{lg }\eta,\phi)$ appropriate for revealing the scaling laws. The
result was instructive: for all experiments of the first group (in
chronological order), specifically: Collins, Coles, Hike,\footnote{The data
were obtained by scanning the graphs in the review by Fernholz and Finley
(1996).} (1978); Erm and Joubert (1991); Smith\footnote{The data were
obtained by scanning the graphs in the review by Fernholz and Finley
(1996).} (1994); Naguib\footnote{The data in digital form were provided to
us by Dr.~M.~Hites.} (1992), and Nagib and Hites\footnote{The data in
digital form were provided to us by Dr.~M.~Hites.} (1995); Krogstad and
Antonia (1998), the data outside the viscous sublayer
$(\mbox{lg }\eta > 1.5)$ have the characteristic shape of a broken line,
shown qualitatively in Figure 1b and quantitatively in Figures 2--6.

Thus, the two straight lines forming the broken line that were revealed in the
$\mbox{lg }\eta,\mbox{lg }\phi$ plane have as equations
\begin{equation}\label{eq11}
{\mbox{(I)}} \ \ \phi=A\eta^\a \ ;
\qquad {\mbox{(II)}} \ \ \phi=B\eta^\b \ . \end{equation}
The coefficients
$A,\a,B,\b$ were obtained by us through statistical processing.

We assume as before that the effective Reynolds number Re has the form
(\ref{eq9}): Re = $U\L/\nu$, where $U$ is the free stream velocity and $\L$
is a length scale. The basic question is, whether one can find in each case
a length scale $\L$ which plays the same role for the intermediate region
(I) of the boundary layer as the diameter does for pipe flow? In other
words, whether it is possible to find a length scale $\L$, perhaps
influenced by individual features of the flow, so that the scaling law
(\ref{eq5}) is valid for the first intermediate region (I)? To answer this
question we have taken the values $A$ and $\a$, obtained by statistical
processing of the experimental data in the first intermediate scaling
region, and then calculated two values $\mbox{ln} \ {\mbox{Re}}_1$,
$\mbox{ln} \ {\mbox{Re}}_2$, by solving the equations suggested by the
scaling law (\ref{eq5}): \begin{equation}\label{eq12}
\frac{1}{\sqrt{3}} \mbox{ln} \ {\mbox{Re}}_1+ \textstyle{\frac 52} =A \ ,
\qquad
\frac{3}{2\mbox{ln} \ {\mbox{Re}}_2}=\a \ . \end{equation}
If these values $\mbox{ln} \ {\mbox{Re}}_1$, $\mbox{ln} \ {\mbox{Re}}_2$
obtained by solving the two different equations (\ref{eq12}) are indeed
close, i.e., if they coincide within experimental accuracy, then the unique
length scale $\L$ can be determined and the experimental scaling law in the
region (I) coincides with the basic scaling law (\ref{eq5}).

Table 1 shows that these values are close, the difference slightly exceeds
3\% in only two cases; in all other cases it is less. Thus, we can
introduce for all these flows the mean Reynolds number
\begin{equation}\label{eq13} {\mbox{Re}}=\sqrt{{\mbox{Re}}_1{\mbox{Re}}_2}
\ , \qquad \mbox{ln} \ \mbox{Re}=\textstyle{\frac 12} (\mbox{ln} \
\mbox{Re}_1+\mbox{ln} \ \mbox{Re}_2) \end{equation}
and consider Re as an estimate of the effective Reynolds number of the
boundary layer flow. Naturally, the ratio Re${}_\th/{\mbox{Re}}=\th/\L$ is
different for different flows.

\newpage
\section{Zero-pressure-gradient boundary layer beneath a turbulent free
stream: The experiments of Hancock and Bradshaw}

The experiments of Hancock
and Bradshaw (1989) revealed a new feature important for our analysis.
Examination of these experimental data suggested that we separate the other
experiments into two groups. In the Hancock and Bradshaw experiments the
free stream was turbulized by a grid in all series, except one. Thus,
processing the data from these experiments we were able not only to compare
the scaling law (\ref{eq5}) with experimental data once again but also to
investigate the influence of the turbulence of the external flow on the
second self-similar intermediate region. The results of the processing are
presented in Table 2 and Figures 7 and 8. In both Table 2 and Figures 7 and
8 the intensity of turbulence is shown by the value of $u'/U$, where
$u'$ is the mean square velocity fluctuation in the free stream.

\begin{center} {\bf Table 2}\\
\begin{tabular}{lrrrrrrrcr}\\
Figure & Re${}_\th$ & $\a$ & A & $\mbox{ln }\mbox{Re}_1$ & $\mbox{ln
}\mbox{Re}_2$ & $\mbox{ln }\mbox{Re}$ & $u'/U$ & Re${}_\th$/Re &
$\b$\bigskip\\ \multicolumn{10}{c} {Hancock, P.E.~and Bradshaw,
P.~(1989)}\\ Fig.8a & 4,680 & 0.140 & 8.66 & 10.67 & 10.71 & 10.69 & 0.0003
& 0.11 & 0.20\\ Fig.8b & 2,980 & 0.138 & 8.77 & 10.86 & 10.91 & 10.88 &
0.024 & 0.06 & 0.18\\ Fig.8c & 5,760 & 0.137 & 8.80 & 10.91 & 10.95 & 10.93
& 0.026 & 0.10 & --\\ Fig.8d & 4,320 & 0.150 & 8.22 & 9.91 & 10.00 & 9.95 &
0.041 & 0.21 & --\\ Fig.8e & 3,710 & 0.122 & 9.49 & 12.11 & 12.30 & 12.20 &
0.040 & 0.02 & --\\ Fig.8f & 3,100 & 0.128 & 9.13 & 11.48 & 11.70 & 11.59 &
0.058 & 0.03 & --\\ Fig.8g & 3,860 & 0.129 & 9.07 & 11.38 & 11.63 & 11.50 &
0.058 & 0.04 & \\ \end{tabular}\end{center}

\medskip
First of all, our processing showed that the first self-similar
intermediate layer is clearly seen in all these experiments, both in the
absence of the external turbulence, and in its presence. The values of
$\mbox{ln} \ {\mbox{Re}}_1$ and $\mbox{ln} \ {\mbox{Re}}_2$ are close. This
means that the basic scaling law
(\ref{eq5}) is valid in the intermediate region
adjacent to the viscous sublayer. At the same
time, the second self-similar region is clearly
observed and well-defined only when the external
turbulence is weak (Figure 8a and to a lesser
extent, Figure 8b) so that the external
turbulence leads to a drastic reduction of the
power $\b$, and even to the deterioration of the
second self-similar intermediate region so that
$\b$ becomes indeterminate. We illustrate the
influence of the free stream turbulence
additionally by Figure 7(b).

The experiments of Hancock and Bradshaw are instructive because they
suggest at least one possible reason for the destruction of the
intermediate self-similar region adjacent to the external flow that is
observed in the experiments of the next group.

\newpage
\section{The remaining group of zero-pressure-gradient boundary layer
experiments}

In this section the results of the processing are presented for all the
remaining series of experiments. For all of them we used the data presented
in the form of graphs in the review of Fernholz and Finley (1996). The
results of the processing are presented in Table 3 and in Figures 9--15.

All the data reveal the self-similar structure in the first intermediate
region adjacent to the viscous sublayer. The scaling laws obtained for this
region give values of $\mbox{ln}$~Re${}_1$ and $\mbox{ln}$~Re${}_2$, close
to each other, although the difference between $\mbox{ln}$~Re${}_1$ and
$\mbox{ln}$~Re${}_2$ is sometimes larger than in the experiments of the
first group. The scaling law (\ref{eq5}) is confirmed by all these
experiments. At the same time, for this group of experiments the second
self-similar structure adjacent to the free stream turns out to be less
clear-cut, if it is there at all. Therefore, for this group of experiments,
we
did not present the estimates
for the values of
$\b$. Nevertheless we note that when it was possible to obtain estimates of
$\b$ they always gave a $\b$ less
than $0.2$. Note also that for all these experiments the number of
experimental points belonging to the region adjacent to the free stream was
less than for the experiments of the first group: this was an additional
argument for our reluctance to show here the second self-similar layer. As
explained in Section 3, we suggest that the turbulence of the external
flow in the experiments of this remaining group was more significant.

\newpage
\section{Checking universality}

The universal form of the scaling law \begin{equation} \label{eq14}
\psi=\frac 1\a \mbox{ln}\left(\frac{2\a\phi}{\sqrt{3}+5\a}\right) =\mbox{ln
}\eta
\end{equation}
gives another way to demonstrate clearly the applicability of the scaling
law (\ref{eq5}) to the first intermediate region of the flow adjacent to
the viscous sublayer. According to relation (\ref{eq14}), in the
coordinates $(\mbox{ln }\eta,\psi)$, all experimental points should
collapse onto the bisectrix of the first quadrant. In Figure 16a are
represented the data of Erm and Joubert (1991), Smith (1994), and Krogstad
and Antonia (1998). It is seen that the data collapse on the bisectrix with
sufficient accuracy to confirm the scaling law (\ref{eq5}). The parameter
$\a$ was calculated according to the formula $\alpha =
(3/2
\mbox{ln}$~Re), $\mbox{ln}$~Re was taken here to be
$(\mbox{ln }{\mbox{Re}}_1+\mbox{ln }{\mbox{Re}}_2)/2$ (see Tables 1 and 3).

In Figure 16b the results of the experiments of Winter and Gaudet (1973),
are presented. These experiments are specially interesting because they
cover a large range of Reynolds numbers (see Table 3). The collapse onto
the bisectrix, although with a larger scatter than for the data presented
in Figure 16a, is clearly demonstrated.

In Figure 16c we present the results of the experiments of Bruns et al
(1992), and Fernholz et al. (1995). Basically they also collapse onto the
bisectrix, although with yet larger scatter, and some systematic deviation
at large $\eta$. This deviation can be explained, at least partially, by
the absence of a sharp outer boundary of the first intermediate region,
unlike the situation in the experiments of the first group.

In Figure 16d are presented the results of all the experiments except those
by Naguib (1992); Naguib and Hites (1995), which will be discussed later,
and the experiments by Winter and Gaudet and Bruns et al. and Fernholz et
al. presented separately in Figures 16b and 16c. As is seen, the
correspondence to the universal form (\ref{eq14}) of the scaling law
(\ref{eq5}) is reasonable. By contrast, Figure 16e,(a) representing the
experiments by Naguib (1992), and Nagib and Hites (1995), shows a
systematic deviation, in fact a parallel shift, from the bisectrix of the
first quadrant. We have already seen such a shift, in the analysis of the
pipe experiments of the Princeton group (Zagarola et al. (1996)); in our
papers on pipe flow (Barenblatt, Chorin, and Prostokishin (1997a, 1997b))
we concluded that the shift was due to the effects of wall roughness, which
increases the effective viscosity. To understand the shift better, we
analyzed also the data in the paper of Krogstad and Antonia (1998) where a
rough wall was used deliberately, albeit for a very large roughness.
Indeed, we found that in these experiments the experimental points lie much
below the bisectrix. Furthermore, in these experiments $\mbox{ln }
\mbox{Re}_1$ and $\mbox{ln } \mbox{Re}_2$ differed significantly, and we
therefore picked the value of $\alpha$ that corresponds to $\mbox{ln }
\mbox{Re}_1$. The result is a pair of lines parallel to the bisectrix but
far below it (Figure 16e,(b)).

More generally, it is very likely that any outside cause that increases the
level of turbulence should also increase the effective viscosity, and thus
shift the points in the $(\mbox{ln } \eta,\psi)$ plane downwards. A case in
point is the set of experiments of Hancock and Bradshaw (1989) discussed
above, where turbulence was created by a grid in the free stream. The
parallel downward shift is indeed observed (Figure 16f), and it is of the
same order of magnitude as the shift in the experiments of Nagib et al.
Note that in the experiments of Nagib et al. the second intermediate region
is intact, and it is therefore likely that the shift in the universal
description of the first intermediate region is due to the disturbance
close to the wall, i.e., to roughness, just as in the experiment of
Zagarola et al. (1996).

\newpage
\section{Conclusion}

The Reynolds-number-dependent scaling law \begin{equation}
\label{eq15}
\phi = \frac{u}{u_*}=
\big(\frac{1}{\sqrt{3}} \mbox{ln} {\mbox{ \ Re}} + {\textstyle{\frac 52}}
\big) \eta^{3/2 \mbox{\footnotesize ln} {\mbox{\footnotesize \ Re}}}
\ ,
\qquad
\eta=\frac{u_*y}{\nu}
\end{equation}
was established earlier for the intermediate region of pipe flows between
the viscous sublayer and the close vicinity of the pipe axis. The Reynolds
number Re was determined as Re $=\bar u \ d/\nu$, where $\bar u$ is the
average velocity, and $d$ the pipe diameter. Attempts (Zagarola and Smits
(1998)) to adjust the constants of the universal logarithmic law so that
this law is valid in a small region of distances from the wall ($y$ less
than $0.07$ of the pipe radius) are immaterial because these data
correspond to the envelope of the family of scaling laws.

In the present work we show that the scaling law (\ref{eq5}) gives an
accurate description of the mean velocity distribution over the
self-similar intermediate region adjacent to the viscous sublayer for a
wide variety of zero-pressure-gradient boundary layer flows. The Reynolds
number is defined as Re $=U\Lambda/\nu$, where $U$ is the free stream
velocity and $\Lambda$ is a length scale which is well defined for all the
flows under investigation.

We also show that under conditions of weak free stream turbulence there
exists a second intermediate self-similar region between the first one,
where the scaling law is valid, and the free stream. This second region
deteriorates under the influence of free stream turbulence.

The validity of the scaling law for boundary layer flows constitutes strong
argument in favor of its validity for a wide class of wall-bounded
turbulent shear flows at large Reynolds numbers. The plotting of the
experimental data in universal coordinates yields a sensitive gauge of the
presence of wall roughness.

Finally,
we feel that the affirmation of the effectiveness of incomplete similarity
and of vanishing-viscosity asymptotics for turbulent shear flows at large
Reynolds numbers has broad implications for other manifestations of
turbulence, e.g.
in jets, wakes, mixing layers, and local structure, and should lead to a
reconsideration of the basic tools used in the study of turbulent flows.

\medskip
\noindent
{\bf Acknowledgements}. The authors would like to thank Professors
P.~Bradshaw and P.-A.~Krogstad for providing them with their recent data.

This work was supported in part by the Applied Mathematics subprogram of
the U.S. Department of Energy under contract DE-AC03-76-SF00098, and in
part by the National Science Foundation under grants DMS94-14631 and
DMS97-32710.

\newpage
\begin{center}
{\bf References}\end{center}

\bigskip\noindent Barenblatt, G. I., 1996. {\it Scaling, Self-Similarity
and Intermediate Asymptotics}, Cambridge University Press.

\bigskip\noindent Barenblatt, G. I. 1991. On the scaling laws (incomplete
self-similarity with respect to Reynolds number) in the developed turbulent
flow in pipes, {\it C.R. Acad. Sci. Paris}, series II, {\bf 313}, 309--312.

\bigskip\noindent Barenblatt, G. I., 1993. Scaling laws for fully developed
shear flows. Part 1: Basic hypotheses and analysis, {\it J. Fluid Mech.}
{\bf 248}, 513--520.

\bigskip\noindent Barenblatt, G. I. and Chorin, A. J., 1996. Small
viscosity asymptotics for the inertial range of local structure and for the
wall region of wall-bounded turbulence, {\it Proc. Nat. Acad. Sciences USA}
{\bf 93}, 6749--6752.

\bigskip\noindent Barenblatt, G. I. and Chorin, A. J., 1997. Scaling laws
and vanishing viscosity limits for wall-bounded shear flows and for local
structure in developed turbulence, {\it Comm. Pure Appl. Math.} {\bf 50},
381--398.

\bigskip\noindent Barenblatt, G. I., Chorin, A. J., Hald, O. H., and
Prostokishin, V. M., 1997. Structure of the zero-pressure-gradient
turbulent boundary layer, 1997. {\it Proc. Nat. Acad. Sciences USA} {\bf
94}, 7817--7819.

\bigskip\noindent Barenblatt, G. I., Chorin, A. J., Prostokishin, V. M.,
1997,a. Scaling laws in fully developed turbulent pipe flow: discussion of
experimental data, {\it Proc. Nat. Acad. Sci. USA} {\bf 94a}, 773--776.

\bigskip\noindent Barenblatt, G. I., Chorin, A. J. and Prostokishin, V. M.,
1997,b. Scaling laws in fully developed turbulent pipe flow, {\it Applied
Mechanics Reviews} {\bf 50}, no.~7, 413--429.

\bigskip\noindent Barenblatt, G. I. and
Prostokishin, V. M., 1993. Scaling laws for fully developed shear flows.
Part 2. Processing of experimental data, {\it J. Fluid Mech.} {\bf 248},
521--529.

\bigskip\noindent
Chorin, A. J., 1988, Scaling laws in the vortex lattice model of
turbulence, {\it Commun. Math. Phys.} {\bf 114}, 167-176.

\bigskip\noindent
Chorin, A. J., 1994. {\it Vorticity and Turbulence}, Springer, New York.

\bigskip\noindent
Chorin, A. J., 1998. New perspectives in turbulence. {\it Quart. Appl.
Math.} {\bf 56}, 767--785.

\bigskip\noindent
Erm, L.P. and Joubert, P. N., 1991. Low Reynolds-number turbulent boundary
layers, {\it J. Fluid Mech.} {\bf 230}, 1--44.

\bigskip\noindent
Fernholz, H. H. and Finley, P. J., 1996. The incompressible
zero-pressure-gradient turbulent boundary layer: an assessment of the data.
{\it Progr. Aerospace Sci.} {\bf 32}, 245--311.

\bigskip\noindent
Hancock, P. E. and Bradshaw, P., 1989. Turbulence structure of a boundary
layer beneath a turbulent free stream, {\it J. Fluid Mech.} {\bf 205},
45--76.

\bigskip\noindent
Krogstad, P.-\AA. and Antonia, R. A., 1998. Surface roughness effects in
turbulent boundary layers, {\it Experiments in Fluids}, in press.

\bigskip\noindent
Landau, L. D. and Lifshits, E. M., 1987. {\it Fluid Mechanics}, Pergamon
Press, New York.

\bigskip\noindent
Monin, A. S. and Yaglom, A. M., 1971. {\it Statistical Fluid Mechanics},
vol. 1, MIT Press, Boston.

\bigskip\noindent
Nagib, H. and Hites, M., 1995. High Reynolds number boundary layer
measurements in the NDF. AIAA paper 95--0786, Reno, Nevada.

\bigskip\noindent
Naguib, A. N., 1992. Inner- and Outer-layer effects on the dynamics of a
turbulent boundary layer, Ph.D. Thesis. Illinois Institute of Technology.

\bigskip\noindent
Nikuradze, J., 1932. Gesetzmaessigkeiten der turbulenten Stroemung in
glatten Rohren, {\it VDI Forschungheft}, No.~356.

\bigskip\noindent
Prandtl, L., 1932. Zur turbulenten Stroemung in Rohren und laengs Platten,
{\it Ergeb. Aerodyn. Versuch.}, Series 4, Goettingen, 18--29.

\bigskip\noindent
Schlichting, H., 1968. {\it Boundary Layer Theory}, McGraw-Hill, New York.

\bigskip\noindent
Spurk, J., 1997. {\it Fluid Mechanics}, Springer, New York.

\bigskip\noindent
von K\'arm\'an, Th., 1930. Mechanische \"Ahnlichkeit und Turbulenz, {\it
Proc. 3rd International Congress for Applied Mechanics}, C.W.Oseen,
W.Weibull, eds. AB Sveriges Litografiska Tryckenier, Stockholm, vol. 1,
pp.85--93.

\bigskip\noindent
Zagarola, M. V., Smits, A. J., Orszag, S. A., and Yakhot, V., 1996.
Experiments in high Reynolds number turbulent pipe flow, AIAA paper
96-0654, Reno, NV.

\bigskip\noindent
Zagarola, M. V. and Smith, A. J., 1998. Mean-flow scaling of turbulent pipe
flow, {\it J. Fluid Mech.} {\bf 373}, 33--79.

\newpage
\begin{center} {\bf Table 1}\\
 \begin{tabular}{lrllrrrcl}\\
 Figure & Re${}_\th$ & $\a$ & A & $\mbox{ln }\mbox{Re}_1$ &
 $\mbox{ln}\mbox{Re}_2$ & $\mbox{ln}\mbox{Re}$ &  Re${}_\th$/Re &
$\b$\bigskip\\
\multicolumn{9}{c} {Collins, D.J., Coles, D.E., and Hiks,
J.W.~(1978)}\\Fig.2(a) & 5,938 & 0.129 & 9.10 & 11.43 & 11.63 & 11.53 &
0.06 & 0.203\\
Fig.2(b) & 6,800 & 0.125 &  9.23 & 11.66 & 12.00 & 11.83 & 0.05 & 0.195\\
Fig.2(c) &  7,880 & 0.123 & 9.41 & 11.97 & 12.21 & 12.09 & 0.04 &
 0.202\bigskip\\
 \multicolumn{9}{c} {Erm, L.P.~and Joubert, P.N.~(1991)}\\  Fig.3(a) & 697
& 0.163 & 7.83 & 9.23 & 9.20 & 9.22 & 0.07 &  0.202\\ Fig.3(b) & 1,003 &
0.159 & 7.96 & 9.46 & 9.43 &  9.45 & 0.08 & 0.192\\ Fig.3(c) & 1,568 &
0.156 & 7.97 &  9.47 & 9.62 & 9.54 & 0.11 & 0.202\\ Fig.3(d) & 2,226 &
 0.148 & 8.26 & 9.98 & 10.14 & 10.06 & 0.10 & 0.214\\  Fig.3(e) & 2,788 &
0.140 & 8.66 & 10.67 & 10.71 & 10.69 &  0.06 & 0.206\bigskip\\
 \multicolumn{9}{c} {Naguib, A.M.~(1992) and Hites, M.~and  Nagib,
H.~(1995)}\\ Fig.4(a) & 4,550 & 0.156 & 7.87 & 9.30  & 9.62 & 9.46 & 0.36 &
0.22\\ Fig.4(b) & 6,240 & 0.148 &  8.24 & 9.94 & 10.14 & 10.04 & 0.27 &
0.20\\ Fig.4(c) &  9,590 & 0.143 & 8.37 & 10.17 & 10.49 & 10.33 & 0.31 &
 0.206\\ Fig.4(d) & 13,800 & 0.131 & 8.94 & 11.15 & 11.45 &  11.30 & 0.17 &
0.193\\ Fig.4(e) & 21,300 & 0.138 & 8.61 &  10.58 & 10.87 & 10.73 & 0.47 &
0.22\\ Fig.4(f) & 29,900 &  0.130 & 8.99 & 11.24 & 11.54 & 11.39 & 0.34 &
0.204\\  Fig.4(g) & 41,800 & 0.124 & 9.30 & 11.78 & 12.10 & 11.94 &  0.27 &
0.201\\ Fig.4(h) & 48,900 & 0.124 & 9.28 & 11.74 &  12.10 & 11.92 & 0.33 &
0.192\bigskip\\
 \multicolumn{9}{c}{Smith, R.W.~(1994)}\\ Fig.5(a) & 4,996  & 0.146 & 8.36
& 10.15 & 10.27 & 10.21 & 0.18 & 0.20\\  Fig.5(b) & 12,990 & 0.129 & 9.19 &
11.59 & 11.63 & 11.61 &  0.12 & 0.167\bigskip\\
 \multicolumn{9}{c} {Krogstad, P.-A.~and Antonia,  R.A.~(1998)}\\ Fig.6 &
12,570 & 0.146 & 8.38 & 10.18 &  10.27 & 10.23 & 0.45 & 0.201
 \end{tabular}
 \end{center}

 \pagebreak

 \begin{center} {\bf Table 3}\\
 \begin{tabular}{lrlrrrrl}\\
 Figure & Re${}_\th$ & $\a$ & A & $\mbox{ln }\mbox{Re}_1$ &  $\mbox{ln
}\mbox{Re}_2$ & $\mbox{ln }\mbox{Re}$ &  Re${}_\th$/Re\bigskip\\
 \multicolumn{8}{c} {Winter, K.G.~and Gaudet, L.~(1973)}\\  Fig.9(a) &
32,150 & 0.133 & 8.86 & 11.02 & 11.32 & 11.17 &  0.45\\ Fig.9(b) & 42,230 &
0.122 & 9.37 & 11.90 & 12.30 &  12.10 & 0.24\\ Fig.9(c) & 77,010 & 0.115 &
10.30 & 13.51 &  13.04 & 13.27 & 0.13\\ Fig.9(d) & 96,280 & 0.107 & 10.56 &
 13.96 & 14.02 & 13.99 & 0.08\\ Fig.9(e) & 136,600 & 0.103 &  10.83 & 14.43
& 14.56 & 14.50 & 0.07\\ Fig.9(f) & 167,600 &  0.101 & 11.20 & 15.07 &
14.85 & 14.96 & 0.05\\ Fig.9(g) &  210,600 & 0.100 & 11.15 & 14.98 & 15.00
& 14.99 &  0.06\bigskip\\
 \multicolumn{8}{c}{ Purtell, L.P., Klebanov, P.S., and  Buckley,
F.T.~(1981)}\\ Fig.10(a) & 1,002 & 0.170 & 7.39 &  8.47 & 8.82 & 8.64 &
0.18\\ Fig.10(b) & 1,837 & 0.164 &  7.62 & 9.14 & 8.87 & 9.00 & 0.23\\
Fig.10(c) & 5,122 &  0.149 & 8.11 & 9.72 & 10.07 & 9.89 & 0.26\bigskip\\
 \multicolumn{8}{c} {Erm, L.P.~(1988)}\\  Fig.11(a) & 2,244 & 0.153 & 8.04
& 9.60 & 9.80 & 9.70 &  0.14\\ Fig.11(b) & 2,777 & 0.154 & 8.13 & 9.75 &
9.74 &  9.75 & 0.16\bigskip\\ \multicolumn{8}{c} {Petrie, H.L.,  Fontaine,
A.A., Sommer, S.T.~and Brungart, T.A.~(1990)}\\  Fig.12 & 35,530 & 0.119 &
9.76 & 12.57 & 12.61 & 12.59 &  0.12\bigskip\\ \multicolumn{8}{c} {Bruns,
J., Dengel, P.,  Fernholz, H.H.~(1992) and}\\
 \multicolumn{8}{c}{Fernholz, H.H., Krause, E., Nockemann,  M. and Schober,
M.~(1995)}\\ Fig.13(a) & 2,573 & 0.151 &  8.46 & 10.32 & 9.93 & 10.13 &
0.10\\ Fig.13(b) & 5,023 &  0.144 & 8.85 & 11.00 & 10.42 & 10.70 & 0.11\\
Fig.13(c) &  7,139 & 0.148 & 8.49 & 10.37 & 10.14 & 10.25 & 0.25\\
 Fig.13(d) & 16,080 & 0.142 & 8.45 & 10.31 & 10.56 & 10.43 &  0.47\\
Fig.13(e) & 20,920 & 0.37 & 8.51 & 10.41 & 10.95 &  10.68 & 0.48\\
Fig.13(f) & 41,260 & 0.132 & 8.63 & 10.62 &  11.36 & 10.98 & 0.70\\
Fig.13(g) & 57,720 & 0.130 & 8.71 &  10.76 & 11.54 & 11.14 & 0.84\bigskip\\
 \multicolumn{8}{c}{Djenidi, L.~and Antonia, R.A.~(1993)}\\  Fig.14(a) &
1,033 & 0.154 & 8.20 & 9.87 & 9.74 & 9.81 &  0.06\\ Fig.14(b) & 1,320 &
0.150 & 8.37 & 10.17 & 10.00 &  10.08 & 0.06\bigskip\\
 \multicolumn{8}{c} {Warnack, D.~(1994)}\\ Fig.15(a) &  2,552 & 0.152 &
8.29 & 10.03 & 9.87 & 9.95 & 0.12\\  Fig.15(b) & 4,736 & 0.149 & 8.20 &
9.87 & 10.07 & 9.97 &  0.22 \end{tabular}\end{center}

\newpage
\noindent
Figure 1.  (a) Schematic representation of the experimental
data in traditional coordinates $\mbox{ln } \eta,\phi$.  (b)
Schematic representation of the experimental data in
$(\mbox{ln } \eta,\mbox{ln } \phi)$ coordinates for
experiments of the first group.

\bigskip
\bigskip

\noindent
Figure 2.  (a) The experiments by Collins, Coles, Hike,
(1978).  $\mbox{Re}_{\theta} = 5,938$.  Both self-similar
intermediate regions (I) and (II) are clearly seen.

\bigskip
\bigskip

\noindent
Figure 2.  (b) The experiments by Collins, Coles, Hike,
(1978).  $\mbox{Re}_{\theta} = 6,800$.  Both self-similar
intermediate regions (I) and (II) are clearly seen.

\bigskip
\bigskip

\noindent
Figure 2.  (c) The experiments by Collins, Coles, Hike,
(1978).  $\mbox{Re}_{\theta} = 7,880$.  Both self-similar
intermediate regions (I) and (II) are clearly seen.

\bigskip
\bigskip

\noindent
Figure 3.  (a) The experiments by Erm and Joubert, (1991).
$\mbox{Re}_{\theta} = 697$.  Both self-similar intermediate
regions (I) and (II) are clearly seen.

\bigskip
\bigskip

\noindent
Figure 3.  (b) The experiments by Erm and Joubert, (1991).
$\mbox{Re}_{\theta} = 1,003$.  Both self-similar
intermediate regions (I) and (II) are clearly seen.

\bigskip
\bigskip

\noindent
Figure 3.  (c) The experiments by Erm and Joubert, (1991).
$\mbox{Re}_{\theta} = 1,568$.  Both self-similar
intermediate regions (I) and (II) are clearly seen.

\bigskip
\bigskip

\noindent
Figure 3.  (d) The experiments by Erm and Joubert, (1991).
$\mbox{Re}_{\theta} = 2,226$.  Both self-similar
intermediate regions (I) and (II) are clearly seen.

\bigskip
\bigskip

\noindent
Figure 3.  (e) The experiments by Erm and Joubert, (1991).
$\mbox{Re}_{\theta} = 2,788$.  Both self-similar
intermediate regions (I) and (II) are clearly seen.

\bigskip
\bigskip

\noindent
Figure 4.  (a) The experiments by Naguib, (1992).
$\mbox{Re}_{\theta} = 4,550$.  Both self-similar
intermediate regions (I) and (II) are clearly seen.

\bigskip
\bigskip

\noindent
Figure 4.  (b) The experiments by Naguib, (1992).
$\mbox{Re}_{\theta} = 6,240$.  Both self-similar
intermediate regions (I) and (II) are clearly seen.

\bigskip
\bigskip

\noindent
Figure 4.  (c) The experiments by Nagib and Hites, (1995).
$\mbox{Re}_{\theta} = 9,590$.  Both self-similar
intermediate regions (I) and (II) are clearly seen.

\bigskip
\bigskip

\noindent
Figure 4.  (d) The experiments by Nagib and Hites, (1995).
$\mbox{Re}_{\theta} = 13,800$.  Both self-similar
intermediate regions (I) and (II) are clearly seen.

\bigskip
\bigskip

\noindent
Figure 4.  (e) The experiments by Nagib and Hites, (1995).
$\mbox{Re}_{\theta} = 21,300$.  Both self-similar
intermediate regions (I) and (II) are clearly seen.

\bigskip
\bigskip

\noindent
Figure 4.  (f) The experiments by Nagib and Hites, (1995).
$\mbox{Re}_{\theta} = 29,900$.  Both self-similar
intermediate regions (I) and (II) are clearly seen.

\bigskip
\bigskip

\noindent
Figure 4.  (g) The experiments by Nagib and Hites, (1995).
$\mbox{Re}_{\theta} = 41,800$.  Both self-similar
intermediate regions (I) and (II) are clearly seen.

\bigskip
\bigskip

\noindent
Figure 4.  (h) The experiments by Nagib and Hites, (1995).
$\mbox{Re}_{\theta} = 48,900$.  Both self-similar
intermediate regions (I) and (II) are clearly seen.

\bigskip
\bigskip

\noindent
Figure 5.  (a) The experiments of Smith, (1994).
$\mbox{Re}_{\theta} = 4,996$.  The first self-similar
intermediate region (I) is clearly seen, the second region
(II) can be revealed.

\bigskip
\bigskip

\noindent
Figure 5.  (b) The experiments of Smith, (1994).
$\mbox{Re}_{\theta} = 12,990$.  The first self-similar
intermediate region (I) is clearly seen, the second region
(II) can be revealed.

\bigskip
\bigskip

\noindent
Figure 6.  The experiments of Krogstad and Antonia, (1998).
$\mbox{Re}_{\theta} = 12,570$.  Both self-similar
intermediate regions (I) and (II) are clearly seen.

\bigskip
\bigskip

\noindent
Figure 7.  (a) The experiments by Hancock and Bradshaw,
(1989) -- a general view.  $\circ$, see Figure 8(a); $+$,
see Figure 8(b); $\times$, see Figure 8(c); $\Box$, see
Figure 8(d); $\Delta$, see Figure 8(e); $\diamond$, see
Figure 8(f); $*$, see Figure 8(g).

\bigskip
\bigskip

\noindent
Figure 7.  (b) The same data as in Figure 7(a) with the
coordinates \newline $\left( x = \mbox{ln } \eta - \frac
{2}{3}
\mbox{ ln } \mbox{Re} \left( \mbox{ln } \frac {\phi}{\frac
{5}{2} + \mbox{\footnotesize ln}\ \mbox{\footnotesize
Re}/\sqrt{3}} \right),\phi\right)$.  The deviations from
the axis $x = 0$ reflect the influence of the turbulence of
the free stream.

\bigskip
\bigskip

\noindent
Figure 8.  (a) The experiments by Hancock and Bradshaw,
(1989).  $\mbox{Re}_{\theta} = 4,680$, $u'/U = 0.0003$.
Both self-similar intermediate regions are clearly seen.

\bigskip
\bigskip

\noindent
Figure 8.  (b) The experiments by Hancock and Bradshaw,
(1989).  $\mbox{Re}_{\theta} = 2,980$, $u'/U = 0.024$.  Both
self-similar intermediate structures (I) and (II) are
clearly seen.

\bigskip
\bigskip

\noindent
Figure 8.  (c) The experiments by Hancock and Bradshaw,
(1989).  $\mbox{Re}_{\theta} = 5,760$, $u'/U = 0.026$.  The
first self-similar intermediate region (I) is clearly seen,
the second is not revealed.

\bigskip
\bigskip

\noindent
Figure 8.  (d) The experiments by Hancock and Bradshaw,
(1989).  $\mbox{Re}_{\theta} = 4,320$, $u'/U = 0.041$.  The
first self-similar intermediate region (I) is clearly seen,
the second is not revealed.

\bigskip
\bigskip

\noindent
Figure 8.  (e) The experiments by Hancock and Bradshaw,
(1989).  $\mbox{Re}_{\theta} = 3,710$, $u'/U = 0.040$.  The
first self-similar intermediate region (I) is clearly seen,
the second is not revealed.

\bigskip
\bigskip

\noindent
Figure 8.  (f) The experiments by Hancock and Bradshaw,
(1989).  $\mbox{Re}_{\theta} = 3,100$, $u'/U = 0.058$.  The
first self-similar intermediate region (I) is seen,
although with a larger scatter.  The second is not revealed.

\bigskip
\bigskip

\noindent
Figure 8.  (g) The experiments of Hancock and Bradshaw,
(1989).  $\mbox{Re}_{\theta} = 3,860$, $u'/U = 0.058$.  The
first self-similar intermediate region (I) is seen,
although with a larger scatter.  The second is not revealed.

\bigskip
\bigskip

\noindent
Figure 9.  (a) The experiments of Winter and Gaudet,
(1973).
$\mbox{Re}_{\theta} = 32,150$.  The first self-similar
intermediate region (I) is seen, although with a larger
scatter.  The second is not revealed.

\bigskip
\bigskip

\noindent
Figure 9.  (b) The experiments of Winter and Gaudet,
(1973).
$\mbox{Re}_{\theta} = 42,230$.  The first self-similar
intermediate region (I) is seen, although with a larger
scatter.  The second is not clearly revealed.

\bigskip
\bigskip

\noindent
Figure 9.  (c) The experiments of Winter and Gaudet,
(1973).
$\mbox{Re}_{\theta} = 77,010$.  The first self-similar
intermediate region (I) is seen, although with a larger
scatter.  The second is not clearly revealed.

\bigskip
\bigskip

\noindent
Figure 9.  (d) The experiments of Winter and Gaudet,
(1973).
$\mbox{Re}_{\theta} = 96,280$.  The first self-similar
intermediate region (I) is seen, although with a larger
scatter.  The second is not clearly revealed.

\bigskip
\bigskip

\noindent
Figure 9.  (e) The experiments of Winter and Gaudet,
(1973).
$\mbox{Re}_{\theta} = 136,600$.  The first self-similar
intermediate region (I) is seen, although with a larger
scatter.  The number of points is not enough to make a
definite estimate for the second region, but the slope
$\beta$ is less than $0.2$.

\bigskip
\bigskip

\noindent
Figure 9.  (f) The experiments of Winter and Gaudet,
(1973).
$\mbox{Re}_{\theta} = 167,600$.  The first self-similar
intermediate region (I) is seen, although with a larger
scatter.  The second is not clearly revealed.

\bigskip
\bigskip

\noindent
Figure 9.  (g) The experiments of Winter and Gaudet,
(1973).
$\mbox{Re}_{\theta} = 210,600$.  The first self-similar
intermediate region (I) is seen, although with a larger
scatter.  The second is not revealed.

\bigskip
\bigskip

\noindent
Figure 10.  (a) The experiments of Purtell, Klebanov and
Buckley, (1981).  $\mbox{Re}_{\theta} = 1,002$.  The first
self-similar region (I) is revealed in spite of the small
number of points.  The second self-similar region is not clearly
revealed.

\bigskip
\bigskip

\noindent
Figure 10.  (b) The experiments of Purtell, Klebanov and
Buckley, (1981).  $\mbox{Re}_{\theta} = 1,837$.  The first
self-similar region (I) is revealed in spite of the small
number of points.  The second self-similar region is
revealed.

\bigskip
\bigskip

\noindent
Figure 10.  (c) The experiments of Purtell, Klebanov and
Buckley, (1981).  $\mbox{Re}_{\theta} = 5,122$.  The first
self-similar region (I) is revealed.  The second
self-similar region is not clearly revealed.

\bigskip
\bigskip

\noindent
Figure 11.  (a) The experiments of Erm, (1988).
$\mbox{Re}_{\theta} = 2,244$.  The first self-similar
region (I) is revealed.  The second self-similar region is
not revealed.

\bigskip
\bigskip

\noindent
Figure 11.  (b) The experiments of Erm, (1988).
$\mbox{Re}_{\theta} = 2,777$.  The first self-similar
region (I) is revealed.  The second self-similar region is
not revealed.

\bigskip
\bigskip

\noindent
Figure 12.  The experiments of Petrie, Fontaine, Sommer and
Brungart, (1990).  $\mbox{Re}_{\theta} = 35,530$.  The first
self-similar region (I) is revealed.  The second
self-similar region is not revealed.

\bigskip
\bigskip

\noindent
Figure 13.  (a) The experiments of Bruns et al., (1992), and
Fernholz et al., (1995).  $\mbox{Re}_{\theta} = 2,573$.  The
first self-similar region (I) is revealed.  The second
self-similar region is not revealed.

\bigskip
\bigskip

\noindent
Figure 13.  (b) The experiments of Bruns et al., (1992), and
Fernholz et al., (1995).  $\mbox{Re}_{\theta} = 5,023$.  The
first self-similar region (I) is revealed.  The second
self-similar region is not revealed.

\bigskip
\bigskip

\noindent
Figure 13.  (c) The experiments of Bruns et al., (1992), and
Fernholz et al., (1995).  $\mbox{Re}_{\theta} = 7,139$.  The
first self-similar region (I) is revealed.  The second
self-similar region is not revealed.

\bigskip
\bigskip

\noindent
Figure 13.  (d) The experiments of Bruns et al., (1992), and
Fernholz et al., (1995).  $\mbox{Re}_{\theta} = 16,080$.
The first self-similar region (I) is revealed.  The second
self-similar region is not revealed.

\bigskip
\bigskip

\noindent
Figure 13.  (e) The experiments of Bruns et al., (1992), and
Fernholz et al., (1995).  $\mbox{Re}_{\theta} = 20,920$.
The first self-similar region (I) is revealed.  The second
self-similar region is not revealed.

\bigskip
\bigskip

\noindent
Figure 13.  (f) The experiments of Bruns et al., (1992), and
Fernholz et al., (1995).  $\mbox{Re}_{\theta} = 41,260$.
The first self-similar region (I) is revealed.  The second
self-similar region is not revealed.

\bigskip
\bigskip

\noindent
Figure 13.  (g) The experiments of Bruns et al., (1992), and
Fernholz et al., (1995).  $\mbox{Re}_{\theta} = 57,720$.
The first self-similar region (I) is revealed.  The second
self-similar region is not revealed.

\bigskip
\bigskip

\noindent
Figure 14.  (a) The experiments of Djenidi and Antonia,
(1993).  $\mbox{Re}_{\theta} = 1,033$.  The first
self-similar region (I) is revealed although with a larger
scatter.  The second self-similar region can be traced.

\bigskip
\bigskip

\noindent
Figure 14.  (b) The experiments of Djenidi and Antonia,
(1993).  $\mbox{Re}_{\theta} = 1,320$.  The first
self-similar region (I) is revealed although with a larger
scatter.  The second self-similar region is not revealed.

\bigskip
\bigskip

\noindent
Figure 15.  (a) The experiments of Warnack, (1994).
$\mbox{Re}_{\theta} = 2,552$.  The first self-similar
region (I) is clearly seen.  The second self-similar region
can be traced.

\bigskip
\bigskip

\noindent
Figure 15.  (b) The experiments of Warnack, (1994).
$\mbox{Re}_{\theta} = 4,736$.  The first self-similar
region (I) is clearly seen.  The second self-similar region
can be traced.

\bigskip
\bigskip

\noindent
Figure 16.  (a) The experiments of Erm and Joubert, (1991)
(*); Smith, (1994) ($\Box$) and Krogstad and Antonia,
(1998) ($\triangleleft$); and Petrie et al., (1990)
($\triangleright$) collapse on the bisectrix of the
first quadrant in accordance with the universal form
$(14)$ of the scaling law $(5)$.

\bigskip
\bigskip

\noindent
Figure 16.  (b) The data of Winter and Gaudet, (1973),
($\nabla$) collapse on the bisectrix of the first quadrant
in accordance with the universal form $(14)$ of the
scaling law $(5)$.

\bigskip
\bigskip

\noindent
Figure 16.  (c) The data of Bruns et al., (1973), and
Fernholz et al., (1995) ($\star$) basically collapse on
the bisectrix of the first quadrant in accordance with the
universal form $(14)$ of the scaling law $(5)$.

\bigskip
\bigskip

\noindent
Figure 16.  (d) The data of all experiments except of those
by Naguib, (1992), and Nagib and Hites, (1995); Bruns et
al., (1992), and Fernholz et al., (1995), collapse on the
bisectrix of the first quadrant in accordance with the
universal form $(14)$ of the scaling law $(5)$:
$(\circ)$ Collins et al., (1978); $(\triangleright)$
Petrie et al., (1990); $(+)$ Erm, (1988); $(\diamond)$
Putell et al., (1981); $(\star)$ Djenidi and Antonia,
(1993);
$(\times)$ Warnack, (1994); $(\triangleleft)$ Krogstad
and Antonia, (1998); $(\nabla)$ Winter and Gaudet, (1973).

\bigskip
\bigskip

\noindent
Figure 16.  (e) (a) The data of Naguib, (1992), and Nagib
and Hites, (1995), show a systematic deviation from the
bisectrix of the first quadrant.

(b) The data of Krogstad an Antonia, (1998), related to
rough walls.  The experimental points lie much lower than
bisectric.  For the evaluation of $\psi$ the value $\alpha
= 3/2 \mbox{ ln } \mbox{Re}_1$ was taken.

\bigskip
\bigskip

\noindent
Figure 16.  (f) The data of Hancock and Bradshaw, (1989),
show the parallel shift from the bisectrix of the same
order as in the experiments by Nagib et al.  $(\cdot)$
Nagib et al., $(\star)$ Hancock and Bradshaw, $u'/U =
0.0003$; $0.024$; $0.026$; $(\times)$ Hancock and Bradshaw,
$u'/U = 0.040,0.041$, $(\circ)$ Hancock and Bradshaw, $u'/U
= 0.058$.

\end{document}